\documentclass[11pt]{article}
\usepackage{amsmath}
\usepackage{amssymb}
\usepackage{url}

\newcommand\nonu{\nonumber}
\newcommand\sa{\smallskipamount}
\newcommand\ma{\medskipamount}
\newcommand\ba{\bigskipamount}
\newcommand\sLP{\\[\sa]}
\newcommand\sPP{\\[\sa]\indent}
\newcommand\mLP{\\[\ma]}
\newcommand\mPP{\\[\ma]\indent}
\newcommand\bLP{\\[\ba]}
\newcommand\bPP{\\[\ba]\indent}
\newcommand\dstyle\displaystyle
\newcommand\CC{\mathbb{C}}
\newcommand\ZZ{\mathbb{Z}}
\newcommand\tha\theta
\newcommand\si\sigma
\newcommand\half{\frac12}
\newcommand\thalf{\tfrac12}
\newcommand\iy\infty
\newcommand\wt{\widetilde}
\newcommand{\qhyp}[5]{\,\mbox{}_{#1}\phi_{#2}\!\left(
  \genfrac{}{}{0pt}{}{#3}{#4};#5\right)}
\newcommand\RHS{right-hand side}

\usepackage[pstarrows]{pict2e}
\usepackage{tikz}
\newcommand{\tikzcircle}[2][red,fill=red]{\tikz[baseline=-0.5ex]
  \draw[thick,#1,radius=#2] (0,0) circle ;}%
\newcommand{\black}{\tikzcircle[black, fill=black]{2.5pt}}%
\newcommand{\white}{\tikzcircle[black, fill=white]{2.5pt}}%

\numberwithin{equation}{section}
\newtheorem{theorem}{Theorem}[section]

\newtheorem{Definition}[theorem]{Definition}

\newtheorem{Remark}[theorem]{Remark}
\newenvironment{remark}{\begin{Remark}\rm}{\end{Remark}}
\newtheorem{Example}[theorem]{Example}

\begin{document}
\title{Charting the $q$-Askey scheme}
\author{Tom H. Koornwinder,
{\small {\tt thkmath@xs4all.nl}}}
\date{\emph{Dedicated to Jasper Stokman on the occasion of his
fiftieth birthday,}\\
\emph{in admiration and friendship}}
\maketitle
\begin{abstract}
Following Verde-Star, Linear Algebra Appl. 627 (2021), we label families
of orthogonal polynomials in the $q$-Askey scheme together with their
$q$-hypergeometric representations by three sequences $x_k, h_k, g_k$
of Laurent polynomials in $q^k$,
two of degree 1 and one of degree 2, satisfying certain constraints.
This gives rise to a precise classification and parametrization of these
families together with their limit transitions. This is displayed in a graphical
scheme. We also describe the four-manifold structure underlying the scheme.
\end{abstract}
\section{Introduction}
The Askey scheme \cite[p.46]{6}, \cite[p.184]{2} and the $q$-Askey scheme
\cite[p.414]{2} display in a graphical way the families of
($q$-)hypergeometric orthogonal polynomials as they occur as limit cases
of the four-parameter top level families: Wilson and Racah polynomials
for the Askey
scheme and Askey--Wilson and $q$-Racah polynomials for the $q$-Askey scheme.
By each arrow to the next lower level one parameter is lost. The bottom
level families no longer depend on parameters.
Since their introduction these schemes have been of great assistance to
everybody who needs to do work with one or more of the families in the
scheme.

These schemes are also expected and partially proven to exist in other contexts,
parallel to the original schemes or generalizing them. These contexts are:
(i) ($q$-)hypergeometric biorthogonal rational functions \cite{20};
(ii) the nonsymmetric case \cite{13}, \cite{14};
(iii) the $q=-1$ case starting with
the Bannai--Ito polynomials \cite[pp.~271--273]{15}, \cite{16};
(iv) generalized (continuous) orthogonal systems~\cite{18};
(v)~$q$-Askey scheme for root system $\textup{BC}_n$, see, among others,
Stokman \cite{19} and references given there.

Still some questions can be posed about the original schemes which, in the
author's opinion, have not
been answered in a satisfactory way until now:
\begin{enumerate}
\item
Are the schemes complete? For answering this question one first needs a
precise criterium for inclusion of a family in the scheme. This criterium
is usually that the orthogonal polynomials should satisfy a Bochner type
property, i.e., that they are eigenfunctions of a second order
linear differential
or ($q$-)difference operator of certain type. However, earlier classifications
\cite{7}, \cite{3} arrive, in the continuous case, at the Askey--Wilson
polynomials being the most general family satisfying the requirements, but do
not give an exhaustive classification of all such families.
A related question is if all limits or specializations from one level to the
level below are present in the scheme. 
\item
Which families deserve an independent status in the scheme and which ones
are just subfamilies of another family?
Several families in the $q$-Askey scheme can be considered as a subfamily,
obtained by restricting the parameters,
of a family higher up in the scheme. Consider for instance the
continuous dual $q$-Hahn polynomials \cite[\S14.3]{2}
and the continuous $q$-Jacobi polynomials \cite[\S14.10]{2},
both being subfamilies of the Askey--Wilson polynomials.
Other subfamilies obtained
by parameter restriction are not in the scheme, and do not even have a name.
What makes the included subfamilies so particular?
\item
Is there a suitable reparametrization of the top level polynomials in the schemes
such that all families lower in the scheme can be obtained by specialization of
parameters? Many arrows in the schemes correspond to taking a limit to 0
or $\iy$ of rescaled polynomials, involving parameter dependent dilation
or translation of the independent variable. It would be nice to simplify this
and make it more uniform. The author \cite{8}
made an attempt in this direction for the Askey scheme. However, there the
formulas for the reparametrization were quite tedious and not very conceptual.
\end{enumerate}

This paper presents, for the $q$-Askey scheme, one possible way to answer these
three questions in a systematic way. Following the ideas by Vinet \& Zhedanov
\cite{9} and Verde-Star \cite{1} one can try to classify
monic orthogonal polynomials
$u_n$ which not only satisfy the Bochner-type property that they are
eigenfunctions of a second order linear $q$-difference operator $L$, so
$Lu_n=h_nu_n$ with the $h_n$ distinct,
but the $u_n$ should also have an expansion
\begin{equation*}
u_n(x)=\sum_{k=0}^n c_{n,k}\,v_k(x),\qquad
v_k(x)=(x-x_0)(x-x_1)\ldots(x-x_{k-1})\;\;(k\ge1),\quad v_0(x)=1.
\end{equation*}
So the $v_n$ are Newton type polynomials.
Now replace the requirement on $L$ to be a second order $q$-difference operator
by the assumption that it acts on the basis of polynomials $v_n$ as
$Lv_n=h_n v_n + g_n v_{n-1}$.
Then it follows that $c_{n,k}=\dstyle\prod_{j=k}^{n-1}\frac{g_{j+1}}{h_n-h_j}$\,.
Finally replace the orthogonality assumption by the property that $xu_n(x)$
is a linear combination of $u_{n+1}(x)$, $u_n(x)$ and (with nonvanishing
coefficient) $u_{n-1}(x)$.

Verde-Star \cite{1}, whom we will follow in this paper, makes the Ansatz
that $h_k$ and $x_k$ are Laurent polynomials in $q^k$ of degree 1 and that
$g_k$ is a Laurent polynomial in $q^k$ of degree~2. The corresponding
$3+3+5=11$ Laurent coefficients then satisfy one trivial relation
because $g_0=0$ and two further relations implied by the three-term recurrence
for the $u_n$. All families in the $q$-Askey scheme \cite[p.414]{2} are
caught by giving the 11 Laurent coefficients, and hence the $x_k,h_k,g_k$,
suitable values. The only exception is the continuous
$q$-Hermite polynomial \cite[\S14.26]{2}. It does not have
an explicit expansion which fits into our framework. 
Apart from this case our method gives a positive answer to the
question whether the $q$-Askey scheme is complete.

It turns out that almost always the distinction between two families in the
scheme can be read off from their different patterns of vanishing Laurent
coefficients (although this does not distinguish between a discrete family and
its continuous analogue). 
These different patterns correspond with different types of
$q$-hypergeometric representations. This answers, in a sense,
the second question.
Furthermore, if we draw an arrow from one family to another in the case that
a suitable
nonzero Laurent coefficient for the first family
becomes zero for the second family, we recover all arrows in the existing
$q$-Askey scheme and find a few more.

Finally, the question about the reparametrization can be answered by starting
with the 11 Laurent coefficients, reduce them by a number of identifications
to a four-manifold, and distinguish lower dimensional submanifolds by
putting one or more suitable Laurent coefficients to zero.

The contents of this paper are as follows. In Section 2 we describe the
general set-up, following Verde-Star \cite{1}, and we illustrate this
for the case of the Askey--Wilson polynomials. In Section 3, the heart of this
paper, we give the resulting scheme. Section 4 describes the manifold
structure associated with the scheme. Finally Section 5 gives some
further perspectives. There are two Appendices. The first one gives
explicit data for families in the scheme. The second one
gives some limit transitions, which are partially missing in
\cite[p.414]{2}.
\paragraph{Acknowledgement}
I thank Paul Terwilliger and the referees for helpful comments.
\paragraph{Note}
For definition and notation of $q$-shifted factorials and
$q$-hypergeometric series we follow \cite[\S1.2]{21}.
We will only need terminating series:
\begin{equation*}
\qhyp rs{q^{-n},a_2,\ldots,a_r}{b_1,\ldots,b_s}{q,z}:=
\sum_{k=0}^n \frac{(q^{-n};q)_k}{(q;q)_k}\,
\frac{(a_2,\ldots,a_r;q)_k}{(b_1,\ldots,b_s;q)_k}\,
\big((-1)^k q^{\half k(k-1)}\big)^{s-r+1}z^k.
\end{equation*}
Here
$(b_1,\ldots,b_s;q)_k:=(b_1;q)_k\ldots(b_s;q)_k$ with
$(b;q)_k:=(1-b)(1-qb)\ldots(1-q^{k-1}b)$ the
$q$-shifted factorial.

For formulas on orthogonal polynomials in the $q$-Askey scheme we
refer to \cite[Chapter 14]{2}.
\section{Askey--Wilson polynomials and Verde-Star's theorem}
Let $u_n(x)$ be an Askey--Wilson polynomial, normalized such
that it is monic in $x=z+z^{-1}$:
\begin{equation}
u_n(x)=p_n\big(\thalf x;a,b,c,d\,|\,q\big)=
\frac{(ab,ac,ad;q)_n}{a^n\,(q^{n-1}abcd;q)_n}\,
\qhyp43{q^{-n},q^{n-1}abcd,az,az^{-1}}{ab,ac,ad}{q,q}.
\label{7}
\end{equation}
We will write some properties of these polynomials in a conceptual form
which we can next use more generally.

Formula \eqref{7} can be rewritten as
\begin{equation}
u_n(x)=\sum_{k=0}^n c_{n,k}\,v_k(x),
\label{1}
\end{equation}
where
\begin{align}
v_k(x)&=(x-x_0)(x-x_1)\ldots(x-x_{k-1})\;\;(k\ge1),\quad v_0(x)=1,
\label{2}\\
x_k&=aq^k+a^{-1}q^{-k},
\label{3}
\end{align}
and
\begin{align}
&c_{n,k}=\prod_{j=k}^{n-1}\frac{g_{j+1}}{h_n-h_j},
\label{4}\\
&h_k=q^{-n}(1-q^n)(1-abcdq^{n-1}),
\label{5}\\
&g_k=a^{-1}q^{-2k+1}(1-abq^{k-1})(1-acq^{k-1})(1-adq^{k-1})(1-q^k).
\label{6}
\end{align}
Note that \eqref{1} expands the Askey--Wilson polynomial in terms of Newton
type polynomials \eqref{2} with nodes \eqref{3}. The expansion coefficients
\eqref{4} are expressed in terms of sequences $h_k$ and $g_k$
given by \eqref{5} and \eqref{6}.
Since we will not consider orthogonality, the only constraints to be imposed
on $q,a,b,c,d\in\CC$ are
\begin{equation*}
q\ne0,\quad 1\notin q^\ZZ,\quad a\ne0,\quad 1\notin abcdq^{\ZZ_{\ge0}}.
\end{equation*}
These constraints make $x_k,h_k,g_k$ well-defined and they let $h_n\ne h_j$ for
$n>0$, $0\le j<n$.

According to \cite[(14.1.7)]{2} there is an explicit second order
$q$-difference operator $L$ such that
\begin{equation}
Lu_n=h_nu_n,\quad n\ge0.
\label{13}
\end{equation}
By \eqref{1} we can also characterize $L$ by its action on the basis
of polynomials $v_n$:
\begin{equation}
Lv_0=h_0v_0,\qquad Lv_n=h_nv_n+g_nv_{n-1},\quad n>0.
\label{12}
\end{equation}

The $h_k$, $x_k$, $g_k$ have the form
\begin{equation}
\label{8}
\begin{split}
&h_k=a_0+a_1 q^k+a_2 q^{-k},\qquad
x_k=b_0+b_1 q^k+b_2 q^{-k},\\
&g_k=d_0+d_1 q^k +d_2 q^{-k}+d_3 q^{2k}+d_4 q^{-2k},\quad
\sum_{i=0}^4 d_i=0.
\end{split}
\end{equation}
Furthermore, we see from \eqref{3}, \eqref{5}, \eqref{6} and \eqref{8} that
\begin{equation}
d_3=q^{-1}a_1b_1,\quad d_4=qa_2b_2.
\label{9}
\end{equation}

Now consider arbitrary sequences
$h_k, x_k, g_k$ ($k\ge0$). Assume $g_0=0$ and
$h_n\ne h_j$ for $n>0$, $0\le j<n$.
Let monic polynomials $v_n$ be given by \eqref{3} and let monic
polynomials $u_n$ of degree~$n$ be expanded in terms of the $v_k$ by
\eqref{1} for certain coefficients $c_{n,k}$. Let $L$ be a linear operator
on the space of polynomials. Then any two of the three formulas \eqref{4},
\eqref{13} and \eqref{12} implies the third formula.

The $q$-case of a recent more general result by Verde-Star
\cite[Theorem 6.1]{1} can be
formulated as follows.
\begin{theorem}
\label{11}
Let $q\ne 0$, $1\notin q^\ZZ$. Let $h_k, x_k, g_k$ have the form \eqref{8}.
Assume that $h_n\ne h_j$ for $n>0$, $0\le j<n$, or equivalently
$a_2\notin a_1 q^{\ZZ_{>0}}$. Let the Newton type polynomials
$v_k$ have the form \eqref{2} ($v_k(x)=(x-b_0)^k$ allowed)
and let the monic polynomials $u_n$ of
degree $n$ be defined by \eqref{1} and \eqref{4}. Then the polynomials
$u_n$ satisfy a three-term recurrence relation
\begin{equation}
xu_n(x)=u_{n+1}(x)+A_n u_n(x)+B_n u_{n-1}(x),\quad n\ge1.
\label{10}
\end{equation}
iff \eqref{9} holds.
\end{theorem}

By \cite[(5.5), (5.6)]{1} the coefficients $A_n$ and $B_n$ in \eqref{10}
are given by
\begin{align}
A_n&=x_n+\frac{g_{n+1}}{h_n-h_{n+1}}-\frac{g_n}{h_{n-1}-h_n}\,,
\label{14}\\
B_n&=\frac{g_n}{h_{n-1}-h_n}\left(\frac{g_{n-1}}{h_{n-2}-h_n}
-\frac{g_n}{h_{n-1}-h_n}+\frac{g_{n+1}}{h_{n-1}-h_{n+1}}+x_n-x_{n-1}\right).
\label{15}
\end{align}
For $n=0$ \eqref{10} and \eqref{14} degenerate to
\begin{equation*}
u_1(x)=x-A_0,\qquad A_0=x_0-\frac{g_1}{h_1-h_0}\,.
\end{equation*}
Verde-Star claims that all families in the $q$-Askey scheme
\cite[Chapter 14]{2}, except for the continuous $q$-Hermite polynomials,
can be obtained in this way. We will make this concrete in the next section.

In Theorem \ref{11} it is allowed that $B_n=0$ for all $n$.
This degenerate case
will certainly happen if $g_n=0$ for all~$n$.
Then $A_n=x_n$ and $u_n=v_n$, clearly
not belonging to a family of orthogonal polynomials. We will not include
this case in our classification.

It is also possible that the $B_n$ are zero because the second factor on
the \RHS\ of \eqref{15} is zero. This case will be included in our
classification.

Finally we may have that $g_n$ vanishes only for some values of $n$.
Let then $n=N+1$ the lowest value of $n$ for which $g_n=0$. Then
$c_{n,k}=0$ if $N<k<n$. If we only consider $u_n$ for $n\le N$ we obtain
one of the finite systems of orthogonal polynomials in the $q$-Askey scheme.

Note that a classification according to Theorem \ref{11} does not use
the usual Bochner type criterium \cite{3} of finding all families of orthogonal
polynomials which are eigenfunctions of a suitable second order $q$-difference
operator. Instead it classifies families of polynomials satisfying a
three-term recurrence relation which have an expansion of specific type
in terms of Newton type polynomials of a specific type. Then there is also
an eigenvalue equation \eqref{12}, involving an operator $L$ defined
by \eqref{13}. For each family it can be shown in an ad hoc way that this
operator $L$ can be written as the second order $q$-difference operator
given in \cite[Chapter 14]{2}. But without having done this computation one
already sees that the obtained numbers $h_n$ are the eigenvalues of $L$ given
in \cite{2}.
\begin{remark}
As sketched in \cite[\S3.3]{9}, if we assume that the $u_n$ satisfy a
three-term recurrence relation \eqref{10} and if we assume \eqref{8} only
for the $h_k$, then \eqref{8} for the $x_k$ and $g_k$ will follow.
\end{remark}
\begin{remark}
The polynomials $u_n$ can be renormalized (under assumptions on the $x_k$)
as polynomials $U_n$ given by \eqref{45}. In this form, and with $g_{N+1}=0$
for some $N$, these polynomials also occur in some of Terwilliger's papers,
in particular, \cite[(10)]{22}, \cite[(85)]{23}. Our $x_i,h_i,g_i$
correspond to Terwilliger's $\tha_i,\tha^*_i,\varphi_i$, respectively.
By \cite[Defs.~7.1, 8.1, 14.1, Theor.~23.2]{23} any Leonard system gives
rise to a three-term recurrence relation, of which renormalized solutions
have the mentioned form. See \cite[\S5]{22} for explicit values of
$\tha_i,\tha^*_i,\varphi_i$.
\end{remark}
\begin{remark}
Geronimus raised the problem to classify orthogonal polynomials $u_n$ and
Newton polynomials $v_k$ which satisfy \eqref{1} with $c_{n,k}=a_{n-k}b_k$.
For an exposition and follow-up of this problem see \cite[\S\S3, 4]{24}
\end{remark}
\section{The $q$-Verde-Star scheme}
Let us again give the data leading to polynomials $u_n$ in the $q$-Askey scheme
according to Theorem~\ref{11}:
\begin{align}
&u_n(x)=\sum_{k=0}^n c_{n,k}\,v_k(x),\quad
v_k(x)=\prod_{j=0}^{k-1}(x-x_j),\quad
c_{n,k}=\prod_{j=k}^{n-1}\frac{g_{j+1}}{h_n-h_j}\,,\label{16}\\
\begin{split}
&x_k=b_2 q^{-k}+b_0+b_1 q^k,\quad
h_k=a_2 q^{-k}+a_0+a_1 q^k,\\
&g_k=d_4 q^{-2k}+d_2 q^{-k}+d_0+d_1 q^k +d_3 q^{2k},
\end{split}\label{27}\\
&\sum_{i=0}^4 d_i=0,\quad d_3=q^{-1}a_1b_1,\quad d_4=qa_2b_2,\label{19}\\
&a_2\ne a_1 q^{\ZZ_{>0}},\quad\mbox{in particular,\quad $a_1$ or $a_2\ne 0$},
\qquad \mbox{$d_i\ne0$ for some $i$}.\label{20}
\end{align}
So everything is determined by the 11 parameters $a_0,a_1,a_2,b_0,b_1,b_2,
d_0,d_1,d_2,d_3,d_4$.
There are several invariances:
\begin{enumerate}
\item
If $a_0\to a_0+\tau$ then $h_k\to h_k+\tau$.
\item
If $a_0,a_1,a_2$ and $d_0,d_1,d_2,d_3,d_4$ are multiplied by $\mu\ne0$ then
$h_k,g_k$ are multiplied by $\mu$.
\item
If $b_0\to b_0+\si$ and $x\to x+\si$ then $x_k\to x_k+\si$
and $u_n(x)\to u_n(x+\sigma)$.
\item
If $b_0,b_1,b_2$ and $d_0,d_1,d_2,d_3,d_4$ are multiplied by $\rho\ne0$ then
$x_k,g_k$ are multiplied by $\rho$, $v_k(x)\to \rho^k v_k(\rho^{-1} x)$ and
$u_n(x)\to \rho^n u_n(\rho^{-1} x)$.
\end{enumerate}
In each case, what is not mentioned remains unchanged.
In items 1 and 2 there is no effect on the $u_n(x)$.
Also the translations and dilations of the independent variable of $u_n$
by items~3 and~4 are not considered as essential changes of a family
of orthogonal polynomials. So the above four items give rise to
four degrees of freedom in the 11 parameters.
Together with the three constraints \eqref{19}
on the parameters, there are
four essential parameters left, in agreement with the number of four parameters
of the Askey--Wilson polynomials.

There are two further remarkable operations which can be performed on the
11 parameters:
\mLP
\emph{$q\leftrightarrow q^{-1}$ exchange}:\quad
$a_1\leftrightarrow a_2$,\;
$b_1\leftrightarrow b_2$,\;
$d_1\leftrightarrow d_2$,\;
$d_3\leftrightarrow d_4$.
\mLP
\emph{$x\leftrightarrow h$ duality}:\quad
$a_0\leftrightarrow b_0$,\;
$a_1\leftrightarrow b_1$,\;
$a_2\leftrightarrow b_2$; assume also that
$b_2\ne b_1 q^{\ZZ_{>0}}$, in particular, $b_1$ or $b_2\ne0$.
This relates $u_n$ given by \eqref{16} to its \emph{dual} $\wt u_n$
given by
\begin{equation}
\wt u_n(x)=\sum_{k=0}^n \wt c_{n,k}\,\wt v_k(x),\quad
\wt v_k(x)=\prod_{j=0}^{k-1}(x-h_j),\quad
\wt c_{n,k}=\prod_{j=k}^{n-1}\frac{g_{j+1}}{x_n-x_j}\,.\label{17}
\end{equation}
If we put
\begin{align}
U_n(x)&=\prod_{j=0}^{n-1} \frac{h_n-h_j}{g_{j+1}}\times u_n(x)=
\sum_{k=0}^n\frac{\prod_{j=0}^{k-1}(h_n-h_j)\times\prod_{j=0}^{k-1}(x-x_j)}
{\prod_{j=1}^k g_j}\,,\label{45}\\
\wt U_m(x)&=\prod_{j=0}^{m-1} \frac{x_m-x_j}{g_{j+1}}\times \wt u_m(x)=
\sum_{k=0}^m\frac{\prod_{j=0}^{k-1}(x_m-x_j)\times\prod_{j=0}^{k-1}(x-h_j)}
{\prod_{j=1}^k g_j}\label{46}
\end{align}
then (see also \cite[(1.9)]{9})
\begin{equation}
U_n(x_m)=\wt U_m(h_n).
\label{47}
\end{equation}

For classification purposes we arrange the 11 parameters in an  array
\begin{equation}
\begingroup\setlength\arraycolsep{2pt}\begin{matrix}&b_2&b_0&b_1&
\\d_4&d_2&d_0&d_1&d_3\\&a_2&a_0&a_1&\end{matrix}\endgroup
\label{18}
\end{equation}
It will turn out that only the vanishing of some of these parameters
determines the families in the scheme. Let
\black\ denote any parameter value (which may be zero)
and \white\ a zero parameter value. So we can represent Askey--Wilson by
\eqref{18} with all entries given by \black\,.
The distribution of \black\ and \white\ in an array \eqref{18} has to
satisfy the following rules:
\begin{enumerate}
\item
If $b_1$ or $a_1$ is \white\ then $d_3$ is \white\,;
if $b_2$ or $a_2$ is \white\ then $d_4$ is \white\\
(because of the second and third formula in \eqref{19}).
\item
$b_0$ and $a_0$ are always \black\\
(because $h_k\to h_k+\tau$ and $x_k\to x_k+\si$ are allowed).
\item
In the second row there are no \white\ ones between two \black\ ones\\
(because it will turn out that only the most left and the most right
nonzero $d_i$ determines the family).
\item
In the third row there are at least two \black\ ones\\
(because of rule 2 and the first part of \eqref{20}).
\item
In the second row there are at least two \black\ ones\\
(because of the first part of \eqref{19} and the second part of \eqref{20}).
\item
Flipping a \black\ into a \white\ causes an arrow between the symbols.\\
(This determines a limit case where a parameter tends to zero. If $b_2$ or
$a_2$ becomes white, then also $d_4$. If $b_1$ or $a_1$ becomes white,
then also $d_3$.)
\item
Reflection with respect to the central column in the black-white array
of form \eqref{18} means $q\leftrightarrow q^{-1}$ exchange.
\item
Reflection with respect to the middle row means
$x\leftrightarrow h$ duality (only possible if
there are at least two \black\ ones in the first row).
\end{enumerate}
In Figure \ref{21} we give half of the scheme according to these rules.
It has to be complemented with the scheme obtained from the present one by
reflecting each diagram with respect to its middle column and preserving all
arrows.
\setlength{\unitlength}{3mm}
\begin{figure}[h]
\centering
\begin{picture}(50,40)
\put(20,36.5)
{\text{$\begingroup\setlength\arraycolsep{2pt}
\begin{matrix}&\black&\black&\black&\\
\black&\black&\black&\black&\black\\&\black&\black&\black&\end{matrix}
\endgroup$}}
\put(20.3,36.1) {\vector(-2,-1){3.2}}
\put(25.0,36.1) {\vector(2,-1){3.2}}
\put(13.4,31.3)
{\text{$\begingroup\setlength\arraycolsep{2pt}
\begin{matrix}&\black&\black&\black&\\
\black&\black&\black&\black&\white\\&\black&\black&\white&\end{matrix}
\endgroup$}}
\put(26.8,31.3)
{\text{$\begingroup\setlength\arraycolsep{2pt}
\begin{matrix}&\black&\black&\white&\\
\black&\black&\black&\black&\white\\&\black&\black&\black&\end{matrix}
\endgroup$}}
\put(17.7,28.7) {\vector(2,-1){3.2}}
\put(27.4,28.7) {\vector(-2,-1){3.2}}
\put(20,24.3)
{\text{$\begingroup\setlength\arraycolsep{2pt}
\begin{matrix}&\black&\black&\white&\\
\black&\black&\black&\black&\white\\&\black&\black&\white&\end{matrix}
\endgroup$}}
\put(16.0,28.8) {\vector(-1,-2){0.9}}
\put(11.8,24.3)
{\text{$\begingroup\setlength\arraycolsep{2pt}
\begin{matrix}&\white&\black&\black&\\
\white&\black&\black&\black&\white\\&\black&\black&\white&\end{matrix}
\endgroup$}}
\put(13.0,30.5) {\vector(-2,-1){5.9}}
\put(3.0,24.3)
{\text{$\begingroup\setlength\arraycolsep{2pt}
\begin{matrix}&\black&\black&\black&\\
\black&\black&\black&\white&\white\\&\black&\black&\white&\end{matrix}
\endgroup$}}
\put(29.5,28.8) {\vector(1,-2){0.9}}
\put(28,24.3)
{\text{$\begingroup\setlength\arraycolsep{2pt}
\begin{matrix}&\black&\black&\white&\\
\black&\black&\black&\white&\white\\&\black&\black&\black&\end{matrix}
\endgroup$}}
\put(32.5,30.5) {\vector(2,-1){5.9}}
\put(37.0,24.3)
{\text{$\begingroup\setlength\arraycolsep{2pt}
\begin{matrix}&\white&\black&\white&\\
\white&\black&\black&\black&\white\\&\black&\black&\black&\end{matrix}
\endgroup$}}
\put(2.4,25) {\vector(-1,-2){2.3}}
\put(-3.0,17)
{\text{$\begingroup\setlength\arraycolsep{2pt}
\begin{matrix}&\black&\black&\black&\\
\black&\black&\white&\white&\white\\&\black&\black&\white&\end{matrix}
\endgroup$}}
\put(6,22.2) {\vector(1,-1){1.8}}
\put(20.6,23) {\vector(-4,-1){11.5}}
\put(4.7,17)
{\text{$\begingroup\setlength\arraycolsep{2pt}
\begin{matrix}&\black&\black&\white&\\
\black&\black&\black&\white&\white\\&\black&\black&\white&\end{matrix}
\endgroup$}}
\put(7.5,22.5) {\vector(2,-1){5}}
\put(14.7,22) {\vector(0,-1){2}}
\put(12.3,17)
{\text{$\begingroup\setlength\arraycolsep{2pt}
\begin{matrix}&\white&\black&\black&\\
\white&\black&\black&\white&\white\\&\black&\black&\white&\end{matrix}
\endgroup$}}
\put(16.2,22.2) {\vector(2,-1){4.6}}
\put(20.0,17)
{\text{$\begingroup\setlength\arraycolsep{2pt}
\begin{matrix}&\white&\black&\black&\\
\white&\white&\black&\black&\white\\&\black&\black&\white&\end{matrix}
\endgroup$}}
\put(16.2,22.6) {\vector(4,-1){11}}
\put(24.5,22) {\vector(2,-1){4}}
\put(37.6,23) {\vector(-2,-1){5.7}}
\put(27.7,17)
{\text{$\begingroup\setlength\arraycolsep{2pt}
\begin{matrix}&\white&\black&\white&\\
\white&\black&\black&\black&\white\\&\black&\black&\white&\end{matrix}
\endgroup$}}
\put(33,22.2) {\vector(2,-1){3.8}}
\put(35.3,17)
{\text{$\begingroup\setlength\arraycolsep{2pt}
\begin{matrix}&\black&\black&\white&\\
\black&\black&\white&\white&\white\\&\black&\black&\black&\end{matrix}
\endgroup$}}
\put(41.5,22.2) {\vector(2,-1){3.8}}
\put(33,22.8) {\vector(4,-1){11}}
\put(43.0,17)
{\text{$\begingroup\setlength\arraycolsep{2pt}
\begin{matrix}&\white&\black&\white&\\
\white&\black&\black&\white&\white\\&\black&\black&\black&\end{matrix}
\endgroup$}}
\put(2,15) {\vector(1,-1){2}}
\put(7.5,14.7) {\vector(-1,-1){2}}
\put(36,15.3) {\vector(-10,-1){28}}
\put(43.6,15.0) {\vector(-10,-1){24}}
\put(3.2,9.7)
{\text{$\begingroup\setlength\arraycolsep{2pt}
\begin{matrix}&\black&\black&\white&\\
\black&\black&\white&\white&\white\\&\black&\black&\white&\end{matrix}
\endgroup$}}
\put(9.8,15) {\vector(3,-1){6.5}}
\put(15,14.6) {\vector(1,-1){2.0}}
\put(28.4,15) {\vector(-4,-1){10.0}}
\put(14.2,9.7)
{\text{$\begingroup\setlength\arraycolsep{2pt}
\begin{matrix}&\white&\black&\white&\\
\white&\black&\black&\white&\white\\&\black&\black&\white&\end{matrix}
\endgroup$}}
\put(23,14.8) {\vector(2,-1){4}}
\put(30,14.9) {\vector(-1,-1){2}}
\put(25.2,9.7)
{\text{$\begingroup\setlength\arraycolsep{2pt}
\begin{matrix}&\white&\black&\white&\\
\white&\white&\black&\black&\white\\&\black&\black&\white&\end{matrix}
\endgroup$}}
\end{picture}
\vspace*{-2cm}
\caption{The $q$-Verde-Star scheme}
\label{21}
\end{figure}
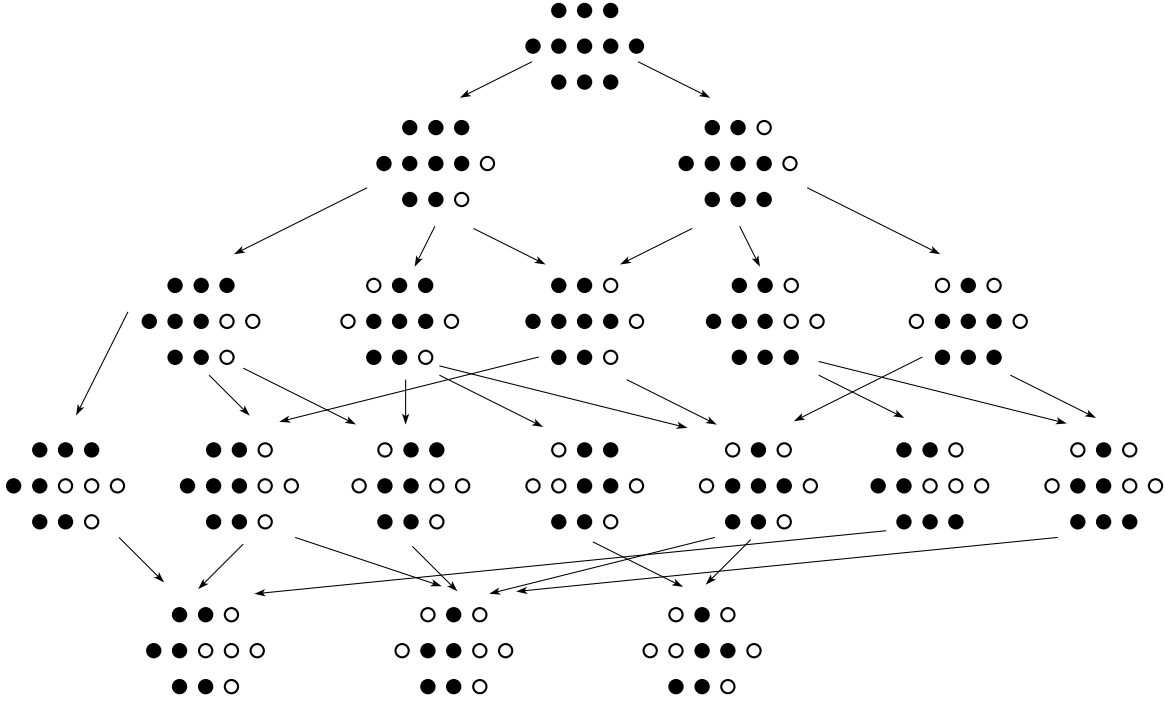

Let us number the rows in the scheme from top to bottom by 1 to 5.
In each row list the successive diagrams from left to right by a, b, $\ldots$\;.
Adding a prime to this notation means a $q\to q^{-1}$ exchange for
the corresponding diagram. For instance \textbf{3c'} is diagram
\textbf{3c} reflected with respect to its central column.
Note that \textbf{1a} = \textbf{1a'} and \textbf{3e} = \textbf{3e'}.
For all other diagrams in Figure \ref{21} the primed counterpart is different
and not in Figure \ref{21}.

Note also the following $x\leftrightarrow h$ dualities:\\
\textbf{1a}, \textbf{3c}, \textbf{4b}, \textbf{5a} are self-dual;\\
\textbf{2a} $\leftrightarrow$ \textbf{2b},
\textbf{3a} $\leftrightarrow$ \textbf{3d},
\textbf{3b} $\leftrightarrow$ \textbf{3b'},
\textbf{4a} $\leftrightarrow$ \textbf{4f},
\textbf{4c} $\leftrightarrow$ \textbf{4d'}
are dual pairs.

The diagrams in Figure \ref{21} correspond with families in 
the $q$-Askey scheme as given in the list below
(numbers given with these families apply to the
corresponding section numbers in \cite[Chapter 14]{2}).
\mLP
\textbf{1a.} Askey--Wilson (1), $q$-Racah (2)
\mLP
\textbf{2a.} continuous dual $q$-Hahn (3), dual $q$-Hahn (7)
\mLP
\textbf{2b.} big $q$-Jacobi (5), $q$-Hahn (6)
\mLP
\textbf{3a.} Al-Salam--Chihara (8), dual $q$-Krawtchouk (17)
\mLP
\textbf{3b.} big $q$-Laguerre (11), $q^{-1}$-Meixner (13),
affine $q$-Krawtchouk (16),
quantum $q^{-1}$-Krawtchouk (14) with
$v_k(x)=\prod_{j=0}^{k-1}(x-b_1 q^j)$.
\mLP
\textbf{3c.} big $q$-Laguerre (11), $q^{-1}$-Meixner (13),
affine $q$-Krawtchouk (16),
quantum $q^{-1}$-Krawtchouk (14) with
$v_k(x)=\prod_{j=0}^{k-1}(x-b_2 q^{-j})$.
\mLP
\textbf{3d.} little $q$-Jacobi (12), $q$-Krawtchouk (15) with
$v_k(x)=\prod_{j=0}^{k-1}(x-b_2 q^{-j})$.
\mLP
\textbf{3e.} little $q$-Jacobi (12), $q$-Krawtchouk (15) with
$v_k(x)=x^k$.
\mLP
\textbf{4a.} continuous big $q$-Hermite (18)
\mLP
\textbf{4b.} $u_n(x)=x^n (bx^{-1};q)_n$,\quad
$v_k(x)=(-1)^k q^{\half k(k-1)}(x;q)_k$.
\mLP
\textbf{4c.} Al-Salam--Carlitz I (24), $q^{-1}$-Al-Salam--Carlitz II (25)
\mLP
\textbf{4d.} little $q$-Laguerre (20), $q^{-1}$-Laguerre (21),
$q^{-1}$-Charlier (23),\quad $v_k(x)=x^k (x^{-1};q)_k$.
\mLP
\textbf{4e.} little $q$-Laguerre (20), $q^{-1}$-Laguerre (21),
$q^{-1}$-Charlier (23),\quad$v_k(x)=x^k$.
\mLP
\textbf{4f.} $q^{-1}$-Bessel (22),\quad
$v_k(x)=(-1)^k q^{\half k(k-1)}(x;q)_k$.
\mLP
\textbf{4g.} $q$-Bessel (22),\quad $v_k(x)=x^k$.
\mLP
\textbf{5a.} $u_n(x)=x^n$,\quad
$v_k(x)=(-1)^k q^{\half k(k-1)}(x;q)_k$.
\mLP
\textbf{5b.} $u_n(x)=x^n (x^{-1};q)_n$,\quad
$v_k(x)=x^k$.
\mLP
\textbf{5c.} $q^{-1}$-Stieltjes--Wigert (27)
\mPP
Figure \ref{21} should be complemented with a similar scheme, where each
diagram is replaced by its primed counterpart and the arrows are preserved.
There are a few arrows from a diagram in the one scheme to a diagram in
the other scheme:
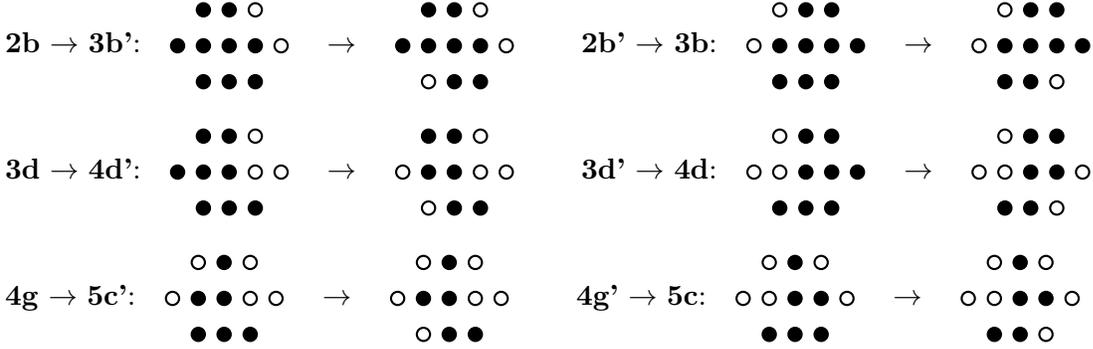
\mLP
\textbf{2b} $\to$ \textbf{3b'}:\quad
$\begingroup\setlength\arraycolsep{2pt}
\begin{matrix}&\black&\black&\white&\\
\black&\black&\black&\black&\white\\&\black&\black&\black&\end{matrix}
\endgroup$
\quad $\to$ \quad
$\begingroup\setlength\arraycolsep{2pt}
\begin{matrix}&\black&\black&\white&\\
\black&\black&\black&\black&\white\\&\white&\black&\black&\end{matrix}
\endgroup$
\qquad
\textbf{2b'} $\to$ \textbf{3b}:\quad
$\begingroup\setlength\arraycolsep{2pt}
\begin{matrix}&\white&\black&\black&\\
\white&\black&\black&\black&\black\\&\black&\black&\black&\end{matrix}
\endgroup$
\quad $\to$ \quad
$\begingroup\setlength\arraycolsep{2pt}
\begin{matrix}&\white&\black&\black&\\
\white&\black&\black&\black&\black\\&\black&\black&\white&\end{matrix}
\endgroup$
\mLP
\textbf{3d} $\to$ \textbf{4d'}:\quad
$\begingroup\setlength\arraycolsep{2pt}
\begin{matrix}&\black&\black&\white&\\
\black&\black&\black&\white&\white\\&\black&\black&\black&\end{matrix}
\endgroup$
\quad $\to$ \quad
$\begingroup\setlength\arraycolsep{2pt}
\begin{matrix}&\black&\black&\white&\\
\white&\black&\black&\white&\white\\&\white&\black&\black&\end{matrix}
\endgroup$
\qquad
\textbf{3d'} $\to$ \textbf{4d}:\quad
$\begingroup\setlength\arraycolsep{2pt}
\begin{matrix}&\white&\black&\black&\\
\white&\white&\black&\black&\black\\&\black&\black&\black&\end{matrix}
\endgroup$
\quad $\to$ \quad
$\begingroup\setlength\arraycolsep{2pt}
\begin{matrix}&\white&\black&\black&\\
\white&\white&\black&\black&\white\\&\black&\black&\white&\end{matrix}
\endgroup$
\mLP
\textbf{4g} $\to$ \textbf{5c'}:\quad
{\text{$\begingroup\setlength\arraycolsep{2pt}
\begin{matrix}&\white&\black&\white&\\
\white&\black&\black&\white&\white\\&\black&\black&\black&\end{matrix}
\endgroup$}}
\quad $\to$ \quad
$\begingroup\setlength\arraycolsep{2pt}
\begin{matrix}&\white&\black&\white&\\
\white&\black&\black&\white&\white\\&\white&\black&\black&\end{matrix}
\endgroup$
\qquad
\textbf{4g'} $\to$ \textbf{5c}:\quad
{\text{$\begingroup\setlength\arraycolsep{2pt}
\begin{matrix}&\white&\black&\white&\\
\white&\white&\black&\black&\white\\&\black&\black&\black&\end{matrix}
\endgroup$}}
\quad $\to$ \quad
$\begingroup\setlength\arraycolsep{2pt}
\begin{matrix}&\white&\black&\white&\\
\white&\white&\black&\black&\white\\&\black&\black&\white&\end{matrix}
\endgroup$
\bLP
\paragraph{Renarks}
\begin{enumerate}
\item
For each diagram in Figure \ref{21} the explicit expressions
of $x_k,h_k,g_k$ for one (continuous) family belonging to this
diagram or its primed counterpart are given in Appendix A.
\item
Since, the first row of a diagram determines the kind
of the Newton type polynomials $v_k$ involved, it can
happen that one family occurs twice in the scheme because it can be expanded
in two different kinds of $v_k$. See \textbf{3b}, \textbf{3c} (big $q$-Laguerre,
affine $q$-Krawtchouk), \textbf{3d},~\textbf{3e} (little $q$-Jacobi,
$q$-Krawtchouk), \textbf{4d}, \textbf{4e} (little $q$-Laguerre,
$q^{-1}$-Laguerre, $q^{-1}$-Charlier), 
\textbf{4f'}, \textbf{4g} ($q$-Bessel).
A family may also have expansions in different $v_k$, which are still of the same
kind. This will not be recognized by our scheme. For instance, with
Askey--Wilson polynomials we may exchange the parameter $a$ with one of the
three other parameters $b,c,d$.
\item
The cases \textbf{4b}, \textbf{5a}, \textbf{5b} are degenerate in the
sense that $u_n$ turns out to be a Newton type polynomial itself which is
expanded in terms of different Newton type polynomials $v_k$.
\item
Most diagrams in the scheme correspond to both a continuous and a discrete
family of orthogonal polynomials.
\item
Figure \ref{21}, when compared with the $q$-Askey scheme on \cite[p.414]{2},
misses some families. The reason is that they are special cases of larger
families, obtained by restriction of parameter values, but such that our
black-white diagrams do not recognize these restrictions. This concerns
(numbers mean again section numbers in \cite[Chapter 14]{2})
continuous $q$-Jacobi (10) as subfamily of Askey--Wilson,
continuous $q$-Laguerre (19)
as subfamily of Al-Salam--Chihara, and discrete $q$-Hermite I, II (28, 29)
as subfamilies of Al-Salam--Carlitz I, II. Similarly, continuous
$q$-Hahn (4) (Askey--Wilson $p_n(x;a,b,c,d\,|\,q)$ with
$a,c$ and $b,d$ pairs of complex conjugates such that $\arg a=\arg c$)
and $q$-Meixner--Pollaczek (9) (Al-Salam--Chihara $Q_n(x;a,b\,|\,q)$ with $a,b$
as a pair of complex conjugates) are not in Figure \ref{21}. (The notation
in \cite[\S\S 14.4, 14.9]{2} for these two classes of polynomials is
confusing.)
\item
The continuous $q$-Hermite polynomials are also missing in our scheme
because their expansion falls
outside the scope of Theorem \ref{11}.
\item
If, in the $q$-Askey scheme on \cite[p.414]{2}, the families mentioned
in the previous two items are omitted, together with the arrows to and from
those families, then all further arrows in that scheme are also present in
our scheme. However, we have some more arrows whcih are missing in
\cite[p.414]{2}. These are (see Appendix B):
\sLP
\textbf{2a} $\to$ \textbf{3b} and \textbf{2a} $\to$ \textbf{3c}:
continuous dual $q$-Hahn $\to$ big $q$-Laguerre,
\sLP
\textbf{3a} $\to$ \textbf{4c}:
Al-Salam--Chihara $\to$ Al-Salam--Carlitz I
\end{enumerate}
\section{The $q$-Verde-Star scheme as a four-manifold}
Fix $q\ne 0$ such that $1\notin q^\ZZ$. We will sketch how the
$q$-Verde-Star scheme can be made into a (complex) four-manifold
having specific submanifolds of lower dimension 3, 2, 1 and 0. We will
ignore the case of a finite system, where $g_{N+1}=0$ for some $N$.

Let us
start with a six-manifold with seven coordinates $a_1,a_2,b_1,b_2,d_0,d_1,d_2$
such that $d_0+d_1+d_2+q^{-1}a_1b_1+qa_2b_2=0$. Also assume that
$a_2\notin a_1 q^{\ZZ_{>0}}$ and $(d_0,d_1,d_2,a_1b_1,a_2b_2)\ne(0,0,0,0,0)$.
Now we make two one-parameter identifications.
Let nothing change if $a_1$, $a_2$, $d_0$, $d_1$, $d_2$ are multiplied by
the same nonzero constant or
if $b_1$, $b_2$, $d_0$, $d_1$, $d_2$ are multiplied by the same nonzero constant.
Then there are several possibilities to put two out of the five coordinates
$a_1$, $a_2$, $b_1$, $b_2$, $d_0$ equal to 1. The three untouched coordinates
among these five, together with $d_1$ or $d_2$ will then provide
local coordinates for our four-manifold. In the generic case all choices are
allowed. However, we can regard the six families in the bottom row of
Figure \ref{21} and its $q\leftrightarrow q^{-1}$ complement as points in
our four-manifold. In a neighbourhood of each of these points we can make
a special choice of four coordinates such that the point has all coordinates
zero and such that any family in the scheme from which that point is
reachable via arrows is a submanifold obtained by putting some of
the coordinates
equal to zero. Below we give the details for the three families in the bottom
row of Figure \ref{21}.
\bPP
\begin{minipage}{5cm}
\qquad\;\;$a_2=b_2=1$
\mLP
$\begin{matrix}
&a_1&b_1&d_0&d_1\\
\textbf{1a}&\black&\black&\black&\black\\
\textbf{2a}&\white&\black&\black&\black\\
\textbf{2b}&\black&\white&\black&\black\\
\textbf{3a}&\white&\black&\black&\white\\
\textbf{3c}&\white&\white&\black&\black\\
\textbf{3d}&\white&\black&\white&\black\\
\textbf{4a}&\white&\black&\white&\white\\
\textbf{4b}&\white&\white&\black&\white\\
\textbf{4f}&\black&\white&\white&\white\\
\textbf{5a}&\white&\white&\white&\white
&&&&\\
&&&&\\
&&&&\\
&&&&
\end{matrix}$
\end{minipage}
\quad
\begin{minipage}{5cm}
\qquad\;\;$a_2=d_0=1$
\mLP
$\begin{matrix}
&a_1&b_1&b_2&d_1\\
\textbf{1a}&\black&\black&\black&\black\\
\textbf{2a}&\white&\black&\black&\black\\
\textbf{2b}&\black&\white&\black&\black\\
\textbf{2b'}&\black&\black&\white&\black\\
\textbf{3a}&\white&\black&\black&\white\\
\textbf{3b}&\white&\black&\white&\black\\
\textbf{3c}&\white&\white&\black&\black\\
\textbf{3e}&\black&\white&\white&\black\\
\textbf{4b}&\white&\white&\black&\white\\
\textbf{4c}&\white&\black&\white&\white\\
\textbf{4e}&\white&\white&\white&\black\\
\textbf{4g}&\black&\white&\white&\white\\
\textbf{5b}&\white&\white&\white&\white\\
\end{matrix}$
\end{minipage}
\quad
\begin{minipage}{5cm}
\qquad\;\;$a_2=d_0=1$
\mLP
$\begin{matrix}
&a_1&b_1&b_2&d_2\\
\textbf{1a}&\black&\black&\black&\black\\
\textbf{2a}&\white&\black&\black&\black\\
\textbf{2b}&\black&\white&\black&\black\\
\textbf{2b'}&\black&\black&\white&\black\\
\textbf{3b}&\white&\black&\white&\black\\
\textbf{3c}&\white&\white&\black&\black\\
\textbf{3e}&\black&\white&\white&\black\\
\textbf{3d'}&\black&\white&\black&\white\\
\textbf{4d}&\white&\black&\white&\white\\
\textbf{4e}&\white&\white&\white&\black\\
\textbf{4g'}&\black&\white&\white&\white\\
\textbf{5e}&\white&\white&\white&\white\\
&&&&
\end{matrix}$
\end{minipage}
\section{Further perspectives}
In a next paper the author will express the coefficients in the relations
defining the Zhedanov algebra associated with a family in the
$q$-Askey scheme (see \cite[(3.2)]{4} with $R=1-\thalf(q+q^{-1})$)
in terms of the 11 parameters
$a_0,a_1,a_2,b_0,b_1,b_2,d_0,d_1,d_2,d_3,d_4$.
It will turn out that vanishing properties of these coefficients are also
a way to distinguish between the families, although the resulting scheme
is slightly different from the scheme in Figure \ref{21}.

Verde-Star \cite{1} introduces polynomials $u_n$ and $v_k$ associated
with sequences $x_k$, $h_k$, $g_k$ as in our \S2, but only assuming that
$x_k$ and $h_k$ are solutions of a certain four-term difference equation
and $g_k$ is a solution of a certain six-term difference equation.
As special cases he has the $q$-case, where $x_k$, $h_k$, $g_k$ have the form
\eqref{8}, the $q=1$ case, and the $q=-1$ case. Earlier, in a somewhat
different approach, these three cases
were examined by Vinet \& Zhedanov \cite{9}.
The author is also planning to write a paper where the $q=1$ case will
be treated systematically and in full detail, just as the $q$-case is
treated in the present paper. We will also deal there with the corresponding
Zhedanov algebra (see \cite[(3.2)]{4} with $R=1$). There will also be need
of a systematic and detailed treatment of the $q=-1$ case. Much material about
this is already available in papers by Vinet \& Zhedanov and coauthors,
see for instance \cite{11}, \cite{10}.

Since the labeling of orthogonal polynomials in the ($q$-)Askey scheme is
by the sequences $x_k$, $h_k$, $g_k$ or by the parameters occurring in
their expansions is so clean, these data may be helpful for recognizing
polynomials in these schemes from the coefficients in the three-term
recurrence relation, assuming that it would be possible to obtain these data
from these
coefficients. See Tcheutia \cite{12} for recent work on this recognition
problem by different methods.

Finally, an approach as in the present paper may be tried in other situations
where (part of) a $q$-Askey scheme occurs, see the examples mentioned
in the second paragraph of the Introduction.
\appendix
\section{Explicit data for the families in Figure \ref{21}}
For each diagram in Figure \ref{21} we give the data of one (continuous)
family belonging to that diagram or its primed counterpart.
Bold numbers like \textbf{1a} follow the convention explained in connection
with Figure \ref{21}.
Numbers in brackets
apply to the corresponding section numbers in \cite[Chapter 14]{2}.
\mLP
\textbf{1a.} Askey--Wilson (1): $u_n(x)=k_n^{-1} p_n(\half x;a,b,c,d\,|\,q)$,
\sPP\,
$x_k=aq^k+a^{-1}q^{-k}$,\quad
$h_k=q^{-k}(1-q^k)(1-abcdq^{k-1})$,
\sPP\,
$g_k=q^{-2k+1}a^{-1}(1-abq^{k-1})(1-acq^{k-1})(1-adq^{k-1})(1-q^k)$,\quad
$k_n=(q^{n-1}abcd;q)_n$.
\mLP
\textbf{2a.} continuous dual $q$-Hahn (3):
$u_n(x)=p_n(\half x;a,b,c\,|\,q)$,
\sPP\,
$x_k=aq^k+a^{-1}q^{-k}$,\quad
$h_k=q^{-k}-1$,\quad
$g_k=q^{-2k+1}a^{-1}(1-abq^{k-1})(1-acq^{k-1})(1-q^k)$.
\mLP
\textbf{2b.} big $q$-Jacobi (5):
$u_n(x)=k_n^{-1} P_n(x;a,b,c;q)$,
\sPP\,
$x_k=q^{-k}$,\quad
$h_k=(1-q^{-k})(-1+q^{k+1}ab)$,\quad
$g_k=q^{1-2k}(1-aq^k)(1-cq^k)(1-q^k)$,
\sPP\,
$k_n=\dstyle\frac{(q^{n+1}ab;q)_n}{(qa;q)_n(qc;q)_n}$\,.
\mLP
\textbf{3a.} Al-Salam--Chihara (8):
$u_n(x)=Q_n(\half x;a,b\,|\,q)$,
\sPP\,
$x_k=aq^k+a^{-1}q^{-k}$,\quad
$h_k=q^{-k}-1$,\quad
$g_k=q^{-2k+1}a^{-1}(1-abq^{k-1})(1-q^k)$.
\mLP
\textbf{3b.} big $q$-Laguerre (11):
$u_n(x)=k_n^{-1} P_n(x;a,b;q)$,\quad
$v_k(x)=x^k (qax^{-1};q)_k$
\sPP\,
$x_k=aq^{k+1}$,\quad
$h_k=q^{-k}-1$,\quad
$g_k=-q^{1-k}b(1-aq^k)(1-q^k)$,\quad
$k_n=\dstyle\frac1{(qa;q)_n(qb;q)_n}$\,.
\mLP
\textbf{3c.} idem, $v_k(x)=(-1)^k q^{-\half k(k-1)} (x;q)_k$
\sPP\,
$x_k=q^{-k}$,\quad
$h_k=q^{-k}-1$,\quad
$g_k=q^{1-2k}(1-aq^k)(1-bq^k)(1-q^k)$.
\mLP
\textbf{3d.} little $q$-Jacobi (12):
$u_n(x)=k_n^{-1} p_n(x;a,b;q)$,\quad
$v_k(x)=(-b)^{-k} q^{-\half k(k+1)} (qbx;q)_k$,
\sPP\,
$x_k=q^{-k-1} b^{-1}$,\quad
$h_k=(1-q^{-k})(-1+q^{k+1}ab)$,\quad
$g_k=(1-q^{-k})(1-b^{-1}q^{-k})$,
\sPP\,
$k_n=(-1)^n q^{-\half n(n-1)}\,\dstyle\frac{(abq^{n+1};q)_n}{(aq;q)_n}$\,.
\mLP
\textbf{3e.} idem, $v_k(x)=x^k$,
\sPP\,
$x_k=0$,\quad
$h_k=(1-q^{-k})(-1+q^{k+1}ab)$,\quad
$g_k=(1-q^{-k})(1-aq^k)$.
\mLP
\textbf{4a.} continuous big $q$-Hermite (18):
$u_n(x)=H_n(\half x;a\,|\,q)$,
\sPP\,
$x_k=aq^k+a^{-1}q^{-k}$,\quad
$h_k=q^{-k}-1$,\quad
$g_k=q^{1-2k}a^{-1}(1-q^k)$.
\mLP
\textbf{4b.} $u_n(x)=x^n (bx^{-1};q)_n$\,,\quad
$x_k=q^{-k}$,\quad
$h_k=q^{-k}-1$,\quad
$g_k=(1-q^{-k})(b-q^{1-k})$.
\mLP
\textbf{4c.} Al-Salam--Carlitz I (24):
$u_n(x)=U_n^{(a)}(x;q)$,
\sPP\,
$x_k=q^k$,\quad
$h_k=q^{-k}-1$,\quad
$g_k=a(1-q^{-k})$.
\mLP
\textbf{4d.} little $q$-Laguerre (20):
$u_n(x)=k_n^{-1} p_n(x;a;q)$,\quad $v_k(x)=x^k (x^{-1};q)_k$,
\sPP\,
$x_k=q^k$,\quad
$h_k=1-q^{-k}$,\quad
$g_k=a(q^k-1)$,\quad
$k_n=\dstyle\frac{(-1)^n q^{-\half n(n-1)}}{(aq;q)_n}$\,.
\sPP\,
Note that $q^{-1}$-Laguerre and little $q$-Laguerre can be essentially
identified with each other
\par\hskip0.2cm by \cite[p.521]{2}.
\mLP
\textbf{4e.} idem, $v_k(x)=x^k$,\quad
$x_k=0$,\quad
$h_k=1-q^{-k}$,\quad
$g_k=q^{-k}(1-aq^k)(1-q^k)$.
\mLP
\textbf{4f'.} $q$-Bessel (22): $u_n(x)=k_n^{-1} y_n(x;a;q)$\quad
$v_k(x)=x^k (x^{-1};q)_k$,
\sPP\,
$x_k=q^k$,\quad
$h_k=(1-q^{-k})(1+aq^k)$,\quad
$g_k=aq^{k-1}(q^k-1)$,\quad
$k_n=(-1)^n q^{-\half n(n-1)}(-aq^n;q)_n$\,.
\mLP
\textbf{4g.} idem, $v_k(x)=x^k$,\quad
$x_k=0$,\quad
$h_k=(1-q^{-k})(1+aq^k)$,\quad
$g_k=q^{-k}-1$.
\mLP
\textbf{5a.} $u_n(x)=x^n$,\quad
$x_k=q^{-k}$,\quad
$h_k=q^{-k}-1$,\quad
$g_k=q^{1-2k}(1-q^k)$.
\mLP
\textbf{5b.} $u_n(x)=k_n^{-1}\,{}_1\phi_0(q^{-n};;q,qx)=x^n (x^{-1};q)_n$\,,
\sPP\,
$x_k=0$,\quad
$h_k=q^{-k}-1$,\quad
$g_k=1-q^{-k}$,\quad
$k_n=(-1)^n q^{-\half n(n-1)}$.
\mLP
\textbf{5c'.} Stieltjes--Wigert (27):
$u_n(x)=k_n^{-1} S_n(x;q)$,
\sPP\,
$x_k=0$,\quad
$h_k=q^k-1$,\quad
$g_k=q^{-k}-1$,\quad
$k_n=\dstyle\frac{(-1)^n q^{n^2}}{(q;q)_n}$\,.
\section{Some explicit limit transitions}
\textbf{2a} $\to$ \textbf{3b}: continuous dual $q$-Hahn $\to$
big $q$-Laguerre (missing in \cite[\S\S14.3, 14.11]{2}).
\begin{equation}
\lim_{a\to0}a^n p_n\big(\thalf a^{-1}x;a,a^{-1}bq,a^{-1}cq\,|\,q\big)
=(bq,cq;q)_n\,P_n(x;b,c;q),
\label{22}
\end{equation}
where continuous dual $q$-Hahn (\cite[(14.3.1)]{2} together with symmetry in
$a,b,c$) and
monic big $q$-Laguerre \cite[(14.11.1)]{2} are respectively represented as
\begin{align}
&p_n\big(\thalf x;a,b,c\,|\,q\big)=
\frac{(ab,bc;q)_n}{b^n}\,
\qhyp32{q^{-n},bz,bz^{-1}}{ab,bc}{q,q},\quad x=z+z^{-1},\quad ab,ac,bc<1,
\label{23}\\*
&(bq,cq;q)_n\,P_n(x;b,c;q)=(-c)^n q^{\half n(n+1)}(bq;q)_n\nonu\\
&\qquad\qquad\qquad\qquad\qquad\qquad
\times\qhyp21{q^{-n},bqx^{-1}}{bq}{q,c^{-1}x},\quad 0<bq<1,\; c<0.
\label{24}
\end{align}
Here and elsewhere in this Appendix,
when we mention conditions on the parameters, these are such that the
coefficient $B_n$ in \eqref{10} is positive, also assuming $A_n$ real.
This assures that the polynomials are orthogonal. 
The conditions above, where the
parameters are assumed real, can be obtained from
\cite[(14.3.5), (14.11.4)]{2}.
By these conditions the passage to the limit in \eqref{22}
can be made while keeping the polynomials orthogonal.
\mLP
\textbf{2a} $\to$ \textbf{3c}:
The same limit \eqref{22} also holds
with other $q$-hypergeometric representations
\cite[(14.3.1), (14.11.1)]{2}:
\begin{align}
&p_n\big(\thalf x;a,b,c\,|\,q\big)=
\frac{(ab,ac;q)_n}{a^n}\,
\qhyp32{q^{-n},az,az^{-1}}{ab,ac}{q,q},\quad x=z+z^{-1},
\label{25}\\
&(bq,cq;q)_n\,P_n(x;b,c;q)=(bq,cq;q)_n\,
\qhyp32{q^{-n},0,x}{bq,cq}{q,q}.
\label{26}
\end{align}
\textbf{3a} $\to$ \textbf{4c}:
Al-Salam--Chihara $\to$ Al-Salam--Carlitz I 
(missing in \cite[\S\S14.8, 14.24]{2}).
\begin{equation}
\lim_{a\to\iy} (a)^{-n} Q_n(\thalf ax;a,ab\,|\,q)=U_n^{(b)}(x;q),
\end{equation}
where Al-Salam--Chihara \cite[(14.8.1)]{2}
and Al-Salam--Carlitz I \cite[(14.24.1)]{2} are respectively represented as
\begin{align}
&Q_n\big(\thalf x;a,b\,|\,q\big)=
\frac{(ab;q)_n}{a^n}\,
\qhyp32{q^{-n},az,az^{-1}}{ab,0}{q,q},\quad x=z+z^{-1},\;ab<1,
\label{40}\\
&U_n^{(b)}(x;q)=(-b)^n q^{\half n(n-1)}\,
\qhyp21{q^{-n},x^{-1}}0{q,qb^{-1}x},\quad b<0.
\end{align}
\textbf{3a} $\to$ \textbf{4b}:
Al-Salam--Chihara $\to$ $x^n(bx^{-1};q)_n$.
\begin{equation}
\lim_{a\to0} a^n Q\big((2a)^{-1}x;a,a^{-1}b\,|\,q\big)=
(b;q)_n \qhyp21{q^{-n},x}b{q,q}=x^n(bx^{-1};q)_n,
\label{41}
\end{equation}
where Al-Salm--Chihara is given by \eqref{40} and the second equality in
\eqref{41} is \cite[(II.6)]{21}.
\mLP
\textbf{2b} $\to$ \textbf{3d}: big $q$-Jacobi $\to$ little $q$-Jacobi
\cite[p.442, Remarks]{2}.
\begin{equation}
\lim_{d\to0}\,(qa)^{-n}\,\frac{(qa;q)_n (-qad;q)_n}{(q^{n+1}ab;q)_n}\,
P_n(x;a,b,1,d;q)=
(-1)^n q^{\half n(n-1)}\,\frac{(qb;q)_n}{(q^{n+1}ab;q)_n}\,p_n(x;b,a;q),
\label{34}
\end{equation}
where big and little $q$-Jacobi are respectively represented by
\cite[p.442, (14.5.1) and Remarks]{2}
\begin{align}
P_n(x;a,b,c,d;q)&=P_n(ac^{-1}qx;a,b,-ac^{-1}d;q)=
\qhyp32{q^{-n},q^{n+1}ab,qac^{-1}x}{qa,-qac^{-1}d}{q,q},\nonu\\
&\qquad\qquad\qquad c,d>0,\;-q^{-1}cd^{-1}<a<q^{-1},\;
-q^{-1}c^{-1}d<b<q^{-1},
\\
p_n(x;a,b;q)&=(-qb)^{-n} q^{-\half n(n-1)}\,\frac{(qb;q)_n}{(qa;q)_n}\,
\qhyp32{q^{-n},q^{n+1}ab,qbx}{qb,0}{q,q},\nonu\\
&\qquad\qquad\qquad\qquad\qquad\qquad\qquad\qquad\qquad\qquad
0<a<q^{-1},\; b<q^{-1}.
\label{36}
\end{align}
\textbf{2b} $\to$ \textbf{3e}:
The same limit \eqref{34} also holds
with other $q$-hypergeometric representations
\cite[(2.39), (2.37)]{5}, \cite[(14.12.1)]{2}:
\begin{align}
P_n(x;a,b,c,d;q)&=\left(-\frac{ad}{bc}\right)^n
\frac{(qb;q)_n (-qbcd^{-1};q)_n}{(qa;q)_n (-qac^{-1}d;q)_n}\,
\qhyp32{q^{-n},q^{n+1}ab,-qbd^{-1}x}{qb,-qbcd^{-1}}{q,q},
\\
p_n(x;a,b;q)&=\qhyp21{q^{-n},q^{n+1}ab}{qa}{q,qx}.
\label{35}
\end{align}
\textbf{3e} $\to$ \textbf{4g}:
little $q$-Jacobi $\to$ $q$-Bessel \cite[(14.12.14)]{2}.
\begin{equation}
\lim_{b\to-\iy} p_n(x;-q^{-1}ab^{-1},b;q)=y_n(x;a;q),
\label{37}
\end{equation}
where the little $q$-Jacobi polynomial is given by \eqref{35}
and the $q$-Bessel function by \cite[(14.22.1)]{2}
\begin{equation}
y_n(x;a;q)=\qhyp21{q^{-n},-aq^n}{0}{q,qx},\qquad a>0.
\label{39}
\end{equation}
\textbf{3d'} $\to$ \textbf{4f'}:
The same limit \eqref{37} also holds with other $q$-hypergeometric
representations
for little $q$-Jacobi and $q$-Bessel:
\begin{align}
p_n(x;a,b;q)&=(-1)^n q^{\half n(n+1)}a^n\,\frac{(bq;q)_n}{(aq;q)_n}\,
\qhyp31{q^{-n},abq^{n+1},x^{-1}}{qb}{q,a^{-1}x},
\label{43}\\
y_n(x;a;q)&=(-1)^n q^{n^2} a^n \qhyp30{q^{-n},-aq^n,x^{-1}}-{q,-a^{-1}x}.
\label{38}
\end{align}
Formula\eqref{43} folllows from \eqref{35} by \cite[(III.8)]{21} and
formula \eqref{38} follows from \eqref{39} by taking the limit $c\to0$ in
\cite[(III.8)]{21}.
\mLP
\textbf{4a} $\to$ \textbf{5a}:
continuous big $q$-Hermite $\to x^n$.
\begin{equation}
\lim_{a\to0}a^{n} H_n\big((2a)^{-1}x;a\,|\,q\big)=
\qhyp21{q^{-n},x}0{q,q}=x^n,
\label{42}
\end{equation}
where continuous big $q$-Hermite is given by \cite[(14.18.1)]{2}
\begin{equation}
H_n\big(\thalf x;a\,|\,q\big)=
a^{-n}\qhyp32{q^{-n},az,az^{-1}}{0,0}{q,q},
\quad x=z+z^{-1},
\end{equation}
and the second equality in \eqref{42} is \cite[(II.6)]{21}.
\mLP
\textbf{4e} $\to$ \textbf{5b}:
little $q$-Laguerre $\to$ $x^n (x^{-1};q)_n$.
\begin{equation}
\lim_{a\to0} (-1)^n q^{\half n(n-1)} (aq;q)_n\,p_n(x;a;q)
=(-1)^n q^{\half n(n-1)}\qhyp10{q^{-n}}-{q,qx}=x^n (x^{-1};q)_n.
\end{equation}
Here little $q$-Laguerre is given by \cite[(14.20.1)]{2}
\begin{equation}
p_n(x;a;q)=\qhyp21{q^{-n},0}{qa}{q,qx}
\label{44}
\end{equation}
and the second equality in \eqref{44} follows from \cite[(II.4)]{21}.

\quad\\
\begin{footnotesize}
\begin{quote}
{ T. H. Koornwinder, Korteweg-de Vries Institute, University of
 Amsterdam,\\
 P.O.\ Box 94248, 1090 GE Amsterdam, The Netherlands;

\vspace{\smallskipamount}
email: }\url{thkmath@xs4all.nl}
\end{quote}
\end{footnotesize}

\end{document}